\documentclass{amsart}
\usepackage{geometry}                
\usepackage{graphicx}
\usepackage{amssymb}
\usepackage{epstopdf,amsmath}
\usepackage{enumitem}

\newtheorem{prop}{Proposition}[section]

\newtheorem{lemma}[prop]{Lemma}
\newtheorem*{lemma*}{Lemma}

\newtheorem{theorem}[prop]{Theorem}
\newtheorem{conj}[prop]{Conjecture}

\newtheorem*{prop*}{Proposition}
\newtheorem*{theorem*}{Theorem}

\newcommand{\WW}{\mathbb{W}}
\newcommand{\PP}{\mathbb{P}}
\newcommand{\FF}{\mathbb{F}}
\newcommand{\CC}{\mathbb{C}}
\newcommand{\TT}{\mathbb{T}}

\newcommand{\RR}{\mathbb{R}}

\newcommand{\SSS}{\mathbb{S}}

\title{Introduction to the proof of the Kakeya conjecture}
\author{Larry Guth}

\begin{document}

\maketitle


In the recent preprint \cite{WZ},  Hong Wang and Joshua Zahl announced a proof of the 3-dimensional Kakeya conjecture.   I have been studying the proof,  and I am going to try to write an introduction to the main ideas in the proof.   (I do not claim to have checked every detail,  but I do think I understand the main ingredients in the argument.)

I hope this note might be helpful for people who are interested in the proof of Kakeya,  and also that it might help a little in the process of digesting and checking the proof.   I am currently working on digesting and checking the proof.   If you are also working on digesting and checking the proof and you would like to be in touch,  feel free to reach out to me.

We will begin by introducing the problem and discussing previous work on it and some of the issues that make it difficult.   Then we will describe the new ideas in the recent work.

\vskip10pt

{\bf Acknowledgements.}  Thanks to Hong Wang and Josh Zahl for many interesting conversations about these ideas,  and also for some helpful comments about this survey.

\section{The Kakeya problem}

Suppose that $\TT$ is a set of $\delta$-tubes in $\RR^n$ with length 1.    We write $U(\TT)$ for $\cup_{T \in \TT} T$.    One version of the Kakeya conjecture in dimension $n$ says 

\begin{conj} \label{conjkakvol} For every $\epsilon >0$, there is a constant $c(n, \epsilon)$ so that if $\TT$ is a set of $\sim \delta^{-(n-1)}$ $\delta$-tubes in $\RR^n$ in $\delta$-separated directions, then 

$$ | U(\TT) | \ge c(n, \epsilon) \delta^\epsilon.$$

\end{conj}

\noindent Wang and Zahl proved this conjecture in dimension $n=3$. 

In fact, they proved a more general estimate called the Kakeya conjecture with convex Wolff axioms.   This therem roughly says that the only way a set of tubes in $\RR^3$ can overlap a lot is by clustering into convex sets.    If $K \subset \RR^3$ is a convex set,  we define

$$ \TT_K := \{ T \in \TT: T \subset K \}. $$

\noindent We define the density of $\TT$ in $K$ as

$$ \Delta(\TT, K) = \frac{ \sum_{T \in \TT_K} |T|}{ |K| }. $$

\noindent The density of $\TT$ in $K$ measures how much the tubes of $\TT$ pack into $K$.   Next we consider the maximum density over all convex sets $K$. 

$$ \Delta_{max}(\TT) := \max_K \Delta(\TT, K). $$

Let us define the typical multiplicity of $\TT$ as

$$ \mu(\TT) = \frac{ \sum_{T \in \TT} |T|}{ | U(\TT) |}. $$

The main theorem of \cite{WZ} is the following version of the Kakeya conjecture which says that $\mu(\TT)$ can only be large when the density of $\TT$ in some convex set is large.

\begin{theorem} \label{main} (Wang-Zahl,  \cite{WZ})  If $\TT$ is a set of $\delta$-tubes in $\RR^3$,   and $\Delta_{max}(\TT) \lessapprox 1$,  then

$$ \mu(\TT) \lessapprox 1. $$

\end{theorem}

This is the Wolff axioms version of the Kakeya conjecture,  which easily implies Conjecture \ref{conjkakvol} for $n=3$.   Given Theorem \ref{main},  a simple sampling argument shows that if $\TT$ is any set of $\delta$-tubes in $\RR^3$,  then $\mu(\TT) \lessapprox \Delta_{max}(\TT)$.   

The proof in \cite{WZ} builds on important contributions by Wolff,  Bourgain,  Katz, Laba, and Tao.   First we will describe these older contributions.   These contributions were interesting and striking,  but a proof of the Kakeya conjecture still looked very difficult.  Second we describe the situation before the recent work and discuss the remaining difficulties.   Then we describe the new contribution in \cite{WZ}.  

The goal of this survey is to describe the general strategy of the proofs,  and we will leave out technical details such as pigeonholing arguments and keeping careful track of small parameters.

\section{A key obstacle: the Heisenberg group}

We begin our story with one of the key obstacles.   There are cousin problems that sound quite similar to the Kakeya conjecture but behave differently.    The Heisenberg group example,  introduced by Katz-Laba-Tao in \cite{KLT},  is an important example of this type.   It concerns the complex analogue of the Kakeya problem with convex Wolff axioms.

We define a complex tube in $\CC^n$ with radius $\delta$ and length $r$  by taking a complex line in $\CC^n$,  intersecting it with a ball of radius $r$,  and then taking its $\delta$-neighborhood.   We define volume on $\CC^n$ by identifying it with $\RR^{2n}$,  so if $T$ is a $\delta$-tube in $\CC^n$ of length 1,  then $|T| \sim \delta^{2(n-1)}$.  The quantities $\Delta(\TT, K)$ and $\Delta_{max}(\TT)$ can be defined as above. 

In $\CC^3$,  there is a set of complex $\delta$-tubes $\TT$ with $\Delta_{max}(\TT) \lessapprox 1$ and $\mu(\TT) \approx \delta^{-1}$.    This example shows that the complex analogue of Theorem \ref{main} is false.  This example is called the Heisenberg group example.   It is based on a quadratic real algebraic hypersurface in $\CC^3$.   There are a couple choices for this hypersurface.  One is the surface $H$ defined by 

$$H = \{ (z_1, z_2, z_3:  |z_1|^2 + |z_2|^2 - |z_3|^2 = 1 \} $$

\noindent  This hypersurface contains many complex lines.   For example,  if $\alpha$ is a unit complex number,  then the line defined by $z_1 = 1$,  $z_3 = \alpha z_2$ lies in $H$.   All these lines pass through the point $(1,0,0) \in H$.   The surface $H$ is very symmetric (it is symmetric under the action of the group $U(2,1)$).   By homogeneity there are infinitely many lines thru every point of $H$.   Taking $\delta$-neighorhoods of these lines gives a set of tubes $\TT$ where $\Delta_{max}(\TT) \lessapprox 1$ and yet $\mu(\TT) \approx \delta^{-1}$.  Ruling out examples like this one will be a key part of our story.

\section{The idea of stickiness}

The proof of the Kakeya conjecture is heavily based on looking at sets of tubes at multiple scales.   Given a scale $\rho \in [\delta, 1]$,  we let $\TT_\rho$ denote the set of $\rho$-tubes formed by thickening the $\delta$-tubes of $\TT$.   When we thicken two distinct $\delta$-tubes,  we may get nearly identical $\rho$-tubes.   When this happens we identify the $\rho$-tubes.   So we typically have $| \TT_\rho | \ll | \TT |$. 

For each $T_\rho \in \TT_\rho$,  we define

$$ \TT_{T_\rho} = \{ T \in \TT: T \subset T_\rho \}. $$

Here is a picture illustrating these different sets of tubes.

\includegraphics[scale=.7]{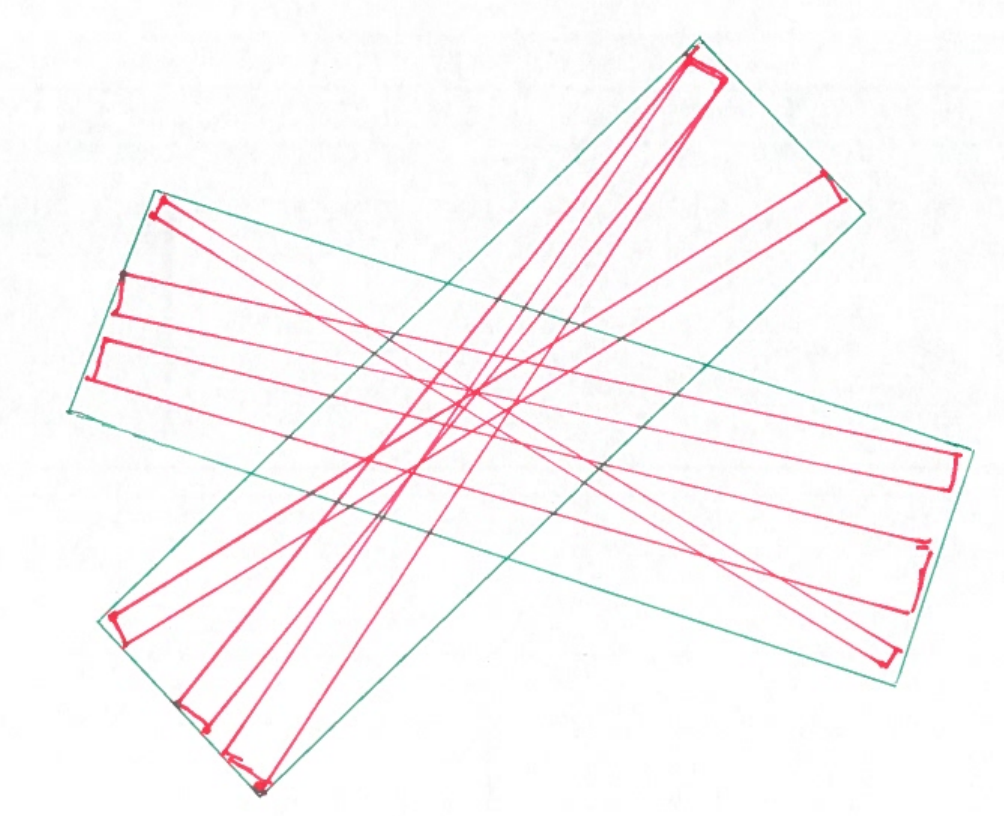}

\noindent The tubes of $\TT$ are in red,  and the tubes of $\TT_\rho$ are in green.   In the picture,  each set $\TT_{T_\rho}$ consists of 3 $\delta$-tubes. 

Our original set of tubes $\TT$ is the disjoint union of $\TT_{T_\rho}$: 

$$ \TT = \bigsqcup_{T_\rho \in \TT_\rho} \TT_{T_\rho}. $$

After some pigeonholing arguments,  we can assume that $|\TT_{T_\rho}|$ is roughly constant for all $T_\rho \in \TT_\rho$.   So for any $T_\rho \in \TT_\rho$,  we get

$$ | \TT | \approx |\TT_{T_\rho}| |\TT_\rho|. $$


By looking at $\TT_{T_\rho}$ and $\TT_\rho$,  we see many different sets of tubes that are all related to our original set of tubes $\TT$.   The key to the proof is to think about all these different sets of tubes and how they relate to each other.

I believe that Tom Wolff was the first person to think about this multi-scale structure in the Kakeya problem,  in unpublished work shortly before his death.   He shared his ideas with Katz, Laba, and Tao, who developed them further in the remarkable paper \cite{KLT}. 

The ``sticky'' case is an important special case of the Kakeya problem.     Because $\Delta_{max}(\TT) \lessapprox 1$,  we see that for every $T_\rho$,  $|\TT_{T_\rho}| \lessapprox \frac{\rho^2}{\delta^2}$.   A Kakeya set is called sticky if this upper bound is roughly tight for all $\rho$ and all $T_\rho \in \TT_\rho$.    In other words,  a Kakeya set is sticky if the tubes of $\TT$ pack into fatter tubes $T_\rho$ as much as they are allowed to do so given that $\Delta_{max}(\TT) \lessapprox 1$.

\newcommand{\gst}{\gamma_{sticky}}

If $\TT$ is near-worst-case among sticky Kakeya sets,  then it must have a some rigid structure.    Recall that $U(\TT)$ denotes $\cup_{T \in \TT} T$.   Suppose that $x \in U(\TT)$.   Typically,  $x$ belongs to $\approx \mu(\TT)$ tubes $T \in \TT$.    Fix a scale $\rho$.   We can now factor this multiplicity $\mu(\TT)$ as

\begin{equation} \label{twoscales}
 \mu (\TT)  \approx \mu(\TT_{T_\rho}) \cdot \# \{T_\rho \in \TT_\rho: x \in U(\TT_{T_\rho}) \} \le \mu(\TT_{T_{\rho}}) \mu(\TT_{\rho}). 
 \end{equation}

Define $\gst$ to be the infimal number so that $\mu(\TT) \lessapprox \delta^{-\gst}$ whenever $\Delta_{max}(\TT) \lessapprox 1$ and $\TT$ is sticky.    Suppose that we are looking at a nearly-worst-case sticky Kakeya set,  so we have $\mu(\TT) \approx \delta^{- \gst}$.    After a linear coordinate transformation,  $\TT_{T_{\rho}}$ can be thought of as a set of $\delta/\rho$ tubes.   Since $\Delta_{max}(\TT) \lessapprox 1$, it follows that $\Delta_{max}(\TT_{T_\rho}) \lessapprox 1$ as well.   Next,  since $\Delta(\TT_{T_\rho},  T_\rho) \approx 1$,  it follows that $\Delta_{max}(\TT_\rho) \lessapprox 1$ as well.      And by the definition of sticky,  $\TT_{T_\rho}$ and $\TT_\rho$ are both sticky also.   So all together we have

$$\delta^{-\gst} \approx \mu(\TT) \lessapprox \mu(\TT_{T_\rho}) \mu(\TT_{\rho}) \lessapprox (\delta/\rho)^{-\gst} \rho^{-\gst} = \delta^{-\gst}. $$

We see that the inequalities above must be approximate equalities.   In particular,  the inequality in (\ref{twoscales}) must be approximately equality.   This gives some striking information!   It tells us that for a typical point $x$ in our Kakeya set, 

$$ \# \{T_\rho \in \TT_\rho: x \in U(\TT_{T_\rho}) \} \approx \mu(\TT_{\rho}).  $$

Consider a point $x \in U(\TT)$ and look at $B(x, \rho)$.   The ball $B(x,\rho)$ lies in $\approx \mu(\TT_\rho) \approx \rho^{-\gst}$ fat tubes $T_\rho \in \TT_\rho$.   For most of these $T_\rho$,  $x \in U(\TT_{T_\rho}) \cap B(x,\rho)$.   That means that the sets $U(\TT_{T_{\rho}}) \cap B(x, \rho)$ must overlap almost perfectly,  as $T_\rho$ varies among the many tubes in $\TT_\rho$ that pass through $B(x,\rho)$. 

Developing these ideas,  Katz, Laba, and Tao showed in \cite{KLT} that such a sticky Kakeya set must have many remarkable structural properties.   For instance,  $U(\TT) \cap B(\delta^{1/2})$ must consist of a union of $\delta \times \delta^{1/2} \times \delta^{1/2}$ plates,  called grains (like the grains in wood).   And all the grains within $B(\delta^{1/2})$ must be parallel to each other,  which is called plaininess. 

The Heisenberg group example in $\CC^3$ is sticky, plainy, and grainy,  as explained in \cite{KLT}.   So we have seen that any sticky Kakeya set must have a great deal of structure in common with the Heisenberg group.   The existence of the Heisenberg group is related to the fact that $\RR$ is a subring of $\CC$,  and the Hausdorff dimension of $\RR$ is half that of $\CC$.   

Katz and Tao developed a strategy to prove that any sticky Kakeya set in $\RR^3$ with large $\mu(\TT)$ must be associated to something like a  fractal approximate subring of $\RR$.   This would be a fractal subset $A \subset \RR$ which is approximately closed under both addition and multiplication.   Bourgain's discretized sum product theorem in \cite{B} shows that no such set exists.   This strategy would then show that any sticky Kakeya set must have $\mu(\TT) \lessapprox 1$.

Katz and Tao did not publish a paper with their strategy,  but Tao put an outline on his blog \cite{T}.   In \cite{WZ1}, Wang and Zahl fleshed out this outline (overcoming significant technical issues),  and proved that every sticky Kakeya set obeys the Kakeya conjecture.

\begin{theorem} (Sticky Kakeya,  \cite{WZ1}) If $\TT$ obeys $\Delta_{max}(\TT) \lessapprox 1$ and if $| \TT_{T_\rho} | \approx (\rho / \delta)^2$ for every $\rho$ and every $T_\rho \in \TT_\rho$,  then

$$ \mu(\TT) \lessapprox 1. $$

\end{theorem}

\noindent  The proof of the sticky Kakeya theorem in \cite{WZ1} uses some important recent developments in projection theory that build on Bourgain's discretized sum-product theorem,  especially \cite{OS1} and \cite{OSW}. 

This work shows a remarkable connection between the sticky case of the Kakeya problem and mathematical structures like the Heisenberg group,  subrings of $\RR$, and sum-product inequalities.   It is certainly interesting mathematics.   But it was not so clear how much progress these results make towards the general Kakeya conjecture.   Is the sticky case a crucial case?  Or is it just a rare special case?

\section{Sticky vs.  not sticky} \label{secstickynot}

Should we expect the ``worst'' Kakeya set to be sticky?  Can we reduce the general problem to the sticky case?  

On the one hand,  it might feel intuitive that the `sticky' case is the worst.   If we are trying to arrange tubes to overlap a lot,  then perhaps we ought to put them close to each other as much as possible.   If we consider a set of tubes $\TT$ without the restriction $\Delta_{max}(\TT) \lessapprox 1$,  then it is easy to arrange that $\mu(\TT)$ is large by packing many tubes into a convex set.   Since the condition $\Delta(\TT,  K) \lessapprox 1$ is limiting us from arranging a lot of overlap,  it may be intuitive that we should pack tubes into convex sets as much as we are allowed to. 

Analysts have considered many cousins of the Kakeya problem.  For many years,  the worse known example for every cousin problem was sticky.   In \cite{B1} Bourgain considered a variation of the Kakeya problem for curved tubes in $\RR^3$.   In this curved version,  the smallest possible volume of $|U(\TT)|$ is $\sim \delta$,  and the worst case example is sticky.   In \cite{KLT},  Katz, Laba, and Tao gave the Heisenberg group example,  which showed that the analogue of the Kakeya problem with convex Wolff axioms is false in $\CC^3$.   The Heisenberg group example is sticky.   In \cite{T2}, Tao noticed that the Kakeya problem with convex Wolff axioms is also false in $\RR^4$.   The worst case example comes from a degree 2 algebraic hypersurface,  and it is also sticky.   It seemed plausible that for a broad class of problems of this type,  the worst case is sticky.

\vskip10pt

However,  it has been difficult to turn this intuition into a rigorous proof.    Katz and Tao and others wondered whether the general Kakeya problem could be reduced to the sticky case,  but they didn't see any way to do it.   In 2017,  in \cite{O},  Orponen proved the sticky case of the Falconer conjecture,  another longstanding problem in geometric measure theory which is a kind of cousin of the Kakeya problem.   This remarkable proof had a lot of influence in the field,  but no one has managed to reduce the general Falconer conjecture to the sticky case. 

In 2019 in \cite{KZ}, Katz and Zahl found a new cousin of the Kakeya problem,  and gave evidence that the worse case example is not sticky.
  They considered the Wolff axioms version of the Kakeya problem over a different ring -- they replaced $\RR$ by the ring $A = \FF_p[x] / (x^2)$.   The ring $A$ has a natural notion of distance with two distinct length scales.   There is a cousin of the Heisenberg group in $A^3$ and it leads to a counterexample to the analogue of Theorem \ref{main}.   But unlike in $\CC^3$,  the Heisenberg group cousin in $A^3$ is {\it not} sticky.   It appears likely that in $A^3$,  the sticky case of (Wolff axiom) Kakeya conjecture is true,  but the general conjecture is false.

As of a couple years ago,  I was quite pessimistic about reducing the general case of Kakeya to the sticky case.  

\vskip10pt

The first indication that the sticky case may play a crucial role in problems of this type was the solution of the Furstenberg set conjecture by Orponen-Shmerkin and Ren-Wang in 2024.   The Furstenberg set conjecture is another important problem in geometric measure theory.   In the late 90s,  Tom Wolff identified the Kakeya problem,  the Falconer problem,  and the Furstenberg set problem as cousin problems involving related issues.
  In 2024,  in \cite{OS2},  Orponen and Shmerkin proved the sticky case of the Furstenberg set conjecture.   Shortly afterwards,  in \cite{RW},  Ren and Wang proved the full Furstenberg conjecture.   Their proof reduces the general problem to the sticky case and a case which is far from sticky,  which they call the semi-well-spaced case.   They resolve the semi-well-spaced case using Fourier analytic methods building on \cite{GSW}.   And they give a short and elegant multiscale argument which reduces the general Furstenberg conjecture to these two cases.

It was quite surprising to me that the general Furstenberg problem reduces to these two cases.  This work gave a hint that the sticky case might be a key case in other problems as well.

\vskip10pt

In 2025,  Wang and Zahl did indeed reduce the general case of the Kakeya problem to the sticky case. 


\section{The $L^2$ method}

Before discussing the work of Wang and Zahl,  let us briefly recall the classical $L^2$ method which plays a supporting role in their work.

If $T_1$ and $T_2$ are two $\delta$-tubes in $\RR^n$ that intersect at a point and the angle between their core lines is $\theta \ge \delta$,  then

\begin{equation} \label{volinttube}
 | T_1 \cap T_2 | \sim \delta^n \theta^{-1}. 
 \end{equation}

If $\TT$ is a set of $\delta$-tubes in $B_1$,  then we can use (\ref{volinttube}) to upper bound

$$ \int_{B_1} | \sum_{T \in \TT} 1_T |^2 = \sum_{T_1, T_2 \in \TT} |T_1 \cap T_2|. $$

Assuming that $| \TT | \approx \delta^{-(n-1)}$ and that $\Delta_{max}(\TT) \lessapprox 1$,  this method gives the sharp bound

$$ \int_{B_1} | \sum_{T \in \TT} 1_T|^2 \lessapprox \delta^{-(n-2)}. $$

\noindent (This bound is sharp when $\TT$ has one tube in each direction and they all go through the origin.)

If $n=2$,  this $L^2$ bound and Cauchy-Schwarz gives a lower bound on $|U(\TT)|$.   If $\TT$ is a set of $\delta$-tubes in $\RR^2$ with $| \TT | \approx \delta^{-1}$ and $\Delta_{max}(\TT) \lessapprox 1$,  this method shows that $|U(\TT)| \gtrapprox 1$,  which is equivalent to $\mu(\TT) \lessapprox 1$.

In higher dimensions,  while this $L^2$ estimate is sharp,  it does not lead to good information about $|U(\TT)|$ or $\mu(\TT)$. 

On the other hand,  this $L^2$ method also works well for slabs in $\RR^3$.   For instance,  using the same method,  we can prove that if $\SSS$ is a set of $\delta \times 1 \times 1$ slabs in $B_1 \subset \RR^3$ with $| \SSS | \sim \delta^{-1}$ and $\Delta_{max}(\SSS) \lessapprox 1$,  then $\mu(\SSS) \lessapprox 1$ and $|U(\SSS)| \gtrapprox 1$.

This method can handle many questions about tubes in $\RR^2$ and slabs in $\RR^3$.   In the proof sketch below,  we will meet a few problems of this type,  and we will mention that they can be handled by the $L^2$ method.

\section{New ideas for the not sticky case}

We now begin to discuss the paper \cite{WZ},  which reduces the general case of the Kakeya conjecture to the sticky case.   We will describe a small variation of the argument from \cite{WZ},  which hopefully simplifies some minor details.   This part of our survey will be more detailed than before, 
but the writeup does not include all details and is not completely precise.   Our goal is to describe the new ideas in the proof.

Define $\beta$ to be the infimal number so that,  whenever $\Delta_{max}(\TT) \lessapprox 1$,  we have

\begin{equation} \label{critexp} \mu(\TT) \lessapprox |\TT|^\beta  \end{equation}

Now suppose that $\TT$ is arranged in a worst case way,  so $\Delta_{max}(\TT) \lessapprox 1$ and $\mu(\TT) \approx | \TT|^\beta$.   Without loss of generality, we can assume that the tubes of $\TT$ lie in the unit ball $B_1$.  Because $\Delta_{max}(\TT) \lessapprox 1$,  we know that $\Delta(\TT,  B_1) \lessapprox 1$ and so $| \TT | \lessapprox \delta^{-2}$.   For these notes,  we focus on the case that $| \TT | \approx \delta^{-2}$.   This case involves most of the ideas in the general proof and is technically simpler.

Our goal is to prove that the worst case is sticky.  Assuming that $\beta > 0$ and that $\TT$ is not sticky, we will prove that $\mu(\TT) \ll | \TT|^\beta$. (We use the notation $\mu(\TT) \ll  |\TT|^\beta$ to mean that $\mu(\TT)$ is much less than $|\TT|^\beta$).  This contradicts our assumption that $\mu(\TT) \approx |\TT|^\beta$, and so we conclude that $\TT$ must be sticky.  But then the sticky Kakeya theorem implies that $\beta = 0$.  

In order to use the definition of $\beta$, we need to relate $\TT$ with other sets of tubes.  We relate $\TT$ to some other set of tubes $\TT'$ with $\Delta_{max}(\TT') \lessapprox 1$ and we use that $\mu(\TT') \lessapprox |\TT'|^\beta$.  Over the course of the proof,  we will bring into play several other sets of tubes $\TT'$.    One key point is to see which sets of tubes $\TT'$ we can bring into play.  

Since $|\TT|\sim \delta^{-2}$,  if $\TT$ is not sticky it means that there is some scale $\rho \in [\delta, 1]$ so that $|\TT_{T_\rho}| \ll (\rho/\delta)^{-2}$ and so $|\TT_\rho| \gg \rho^{-2}$.   We will assume that $\TT$ is not-sticky-at-all-scales,  meaning that for every $\rho$ in the range $\delta \ll \rho \ll 1$,  we have

$$ | \TT_{T_\rho} | \ll (\rho/\delta)^{-2} \textrm{ and } | \TT_\rho | \gg \rho^{-2}. $$

\noindent At the very end,  we will discuss how to reduce the not-sticky case to the not-sticky-at-all-scales case.

\subsection{Some ideas that don't work}

Before turning to the ideas in \cite{WZ},  let me mention some ideas that don't work.   First,  you might imagine trying to slide the tubes of $\TT$ around in order to make $\TT$ more sticky while also increasing $\mu(\TT)$.   I have no idea how to do that.   A hypothetical Kakeya set has a very delicate structure that causes tubes to overlap a lot,  and I don't see any way of sliding the tubes around without breaking the structure and decreasing $\mu(\TT)$.  

Next let's consider the multiscale argument that we used at the start of the sticky case.   Let $\rho \in [\delta, 1]$ be an intermediate scale and consider $\TT_{T_\rho}$ and $\TT_\rho$.   As in (\ref{twoscales}) above,  we can factor the multiplicity as

\begin{equation} \label{twoscales2}
 \mu (\TT)  \approx \mu(\TT_{T_\rho}) \cdot \# \{T_\rho \in \TT_\rho: x \in U(\TT_{T_\rho}) \} \le \mu(\TT_{T_{\rho}}) \mu(\TT_{\rho}). 
 \end{equation}

Because we are not in the not-sticky-at-all-scales case,  we know that for each scale $\rho$ with $\delta \ll \rho \ll 1$, 

$$ | \TT_{T_\rho} | \ll (\rho/ \delta)^{-2}, $$

$$ | \TT_\rho |  \gg \rho^{-2}. $$

\noindent Since $\TT_{T_\rho} \subset \TT$,  we know that $\Delta_{max}(\TT_{T_\rho}) \le \Delta_{max}(\TT) \lessapprox 1$,  and so by the definition of $\beta$,  we can bound

\begin{equation} \label{multfine} \mu(\TT_{T_\rho}) \lessapprox |\TT_{T_{\rho}}|^\beta. \end{equation}

On the other hand,  $\Delta_{max}(\TT_\rho)$ is NOT $\lessapprox 1$,  since $| \TT_\rho| \gg \rho^{-2}$.   This is a key feature of the non-sticky case.   We are thus led to consider sets of tubes with high density.   Given that $\Delta_{max}(\TT) \lessapprox 1$ and $| \TT | \approx \delta^{-2}$,  we get the following information about the density of $\TT_\rho$: for any convex $K \subset B_1$, 

\begin{equation} \label{Trhofrost} \Delta(\TT_\rho, K) \lessapprox \Delta(\TT_\rho, B_1) \sim \rho^2 | \TT_\rho|.  \end{equation}

The condition (\ref{Trhofrost}) plays an important role in the story.   We say that the set of tubes $\TT_\rho$ is Frostman if it obeys (\ref{Trhofrost}).   Wang and Zahl picked this name because the bound in (\ref{Trhofrost}) is similar to the bound that appears in Frostman's lemma in geometric measure theory.

Given that $\beta$ is the best exponent so that $\mu(\TT) \lessapprox |\TT|^\beta$ when $\Delta_{max}(\TT) \lessapprox 1$,  and given that 
$\TT_\rho$ obeys the Frostman condition (\ref{Trhofrost}), 
what can we say about $\mu(\TT_\rho)$ ?  One scenario is that the tubes of $\TT_\rho$ may pack into thick tubes $T_\sigma$,  so that $\Delta_{max}(\TT_\sigma) \approx 1$ and for each $T_\sigma \subset \TT_\sigma$,  $\TT_\rho$ contains every $\rho$-tube in $T_\sigma$.   In this scenario,  after a short computation,  we get

\begin{equation} \label{mufrostman}
\mu(\TT_\rho) \approx (\rho^{-2})^\beta (\rho^2 | \TT_\rho|)^{1 - \beta} 
\end{equation} 

\noindent Wang and Zahl proved that this is the largest possible value of $\mu(\TT_\rho)$:

\begin{lemma} \label{lemmahighdensity} (High density lemma) Suppose that the exponent $\beta$ is as defined above.   If $\TT_\rho$ obeys the Frostman condition  (\ref{Trhofrost}),  then the multiplicty of $\TT_\rho$ is bounded by

$$ \mu(\TT_\rho) \lessapprox (\rho^{-2})^\beta (\rho^2 | \TT_\rho|)^{1 - \beta} . $$

\end{lemma}

The proof of this lemma is based on a subtle multiscale analysis,  and it also uses the sticky Kakeya theorem.   We will return to it in Section \ref{sechighdens}.

Now that we have thought through $\mu(\TT_{T_\rho})$ and $\mu(\TT_\rho)$,  we can return to (\ref{twoscales2}).   Plugging in our bounds,  we see that the right-hand side is much larger than the left-hand side.   In the sticky case,  the two sides were comparable which gave us a lot of rigidity.   In the non-sticky case,  we don't have that rigidity,  and this approach doesn't seem to give really helpful information.

However,  during this discussion,  we did meet the high density lemma,  which will be crucial in the main argument.  


\subsection{Looking at $\TT$ inside a small ball}

Let $\rho$ be an intermediate scale with $\delta \ll \rho \ll 1$.   We will ultimately pick $\rho$ to be quite close to $\delta$.

To bound $\mu(\TT)$ in the sticky case,  we considered  $\TT_{T_\rho}$ and $\TT_\rho$.  In the non-sticky case, Wang and Zahl also consider a new set of tubes formed by intersecting tubes of $\TT$ with a smaller ball $B \subset B_1$.   

To set this up,  let's first let's think about how the tubes of $\TT_{T_\rho}$ intersect each other.   Consider the multiplicity function 

$$\mu_{\TT_{T_\rho}}(x) = \sum_{T \in \TT_{T_\rho}} 1_T(x). $$

\newcommand{\Seg}{\textrm{Seg}}

\noindent Note that $\mu_{\TT_{T_\rho}}$ is morally constant on shorter tubes,  of dimensions $\delta \times \delta \times \delta/\rho$.    Now let $B$ be a ball of radius $\delta /\rho$ and consider how these shorter tubes overlap inside of $B$.  Let $\TT_{T_\rho, B}$ be the set of these shorter tubes in $B$.  So each tube $T_B \in \TT_{T_\rho,  B}$ is a $\delta \times \delta \times \delta/\rho$ tube in $B$ which lies in $\approx \mu(\TT_{T_\rho})$ tubes of $\TT_{T_\rho}$.    Here is a picture.

\includegraphics[scale=.6]{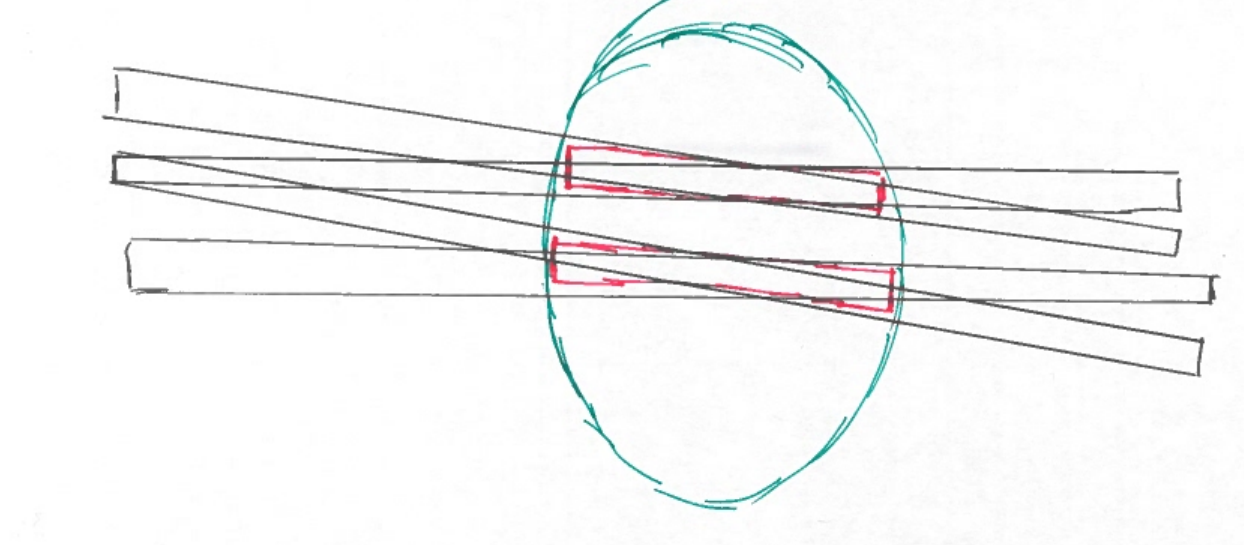}

In this picture,  the long tubes belong to $\TT_{T_\rho}$,  the short red tubes belong to $\TT_{T_\rho, B}$,  and the green disk is $B$. 

Next we define

$$ \TT_B = \bigcup_{T_\rho \in \TT_\rho,  T \cap B \not= \emptyset} \TT_{T_\rho, B}. $$

Here is a picture showing $\TT_B$:

\includegraphics[scale=.5]{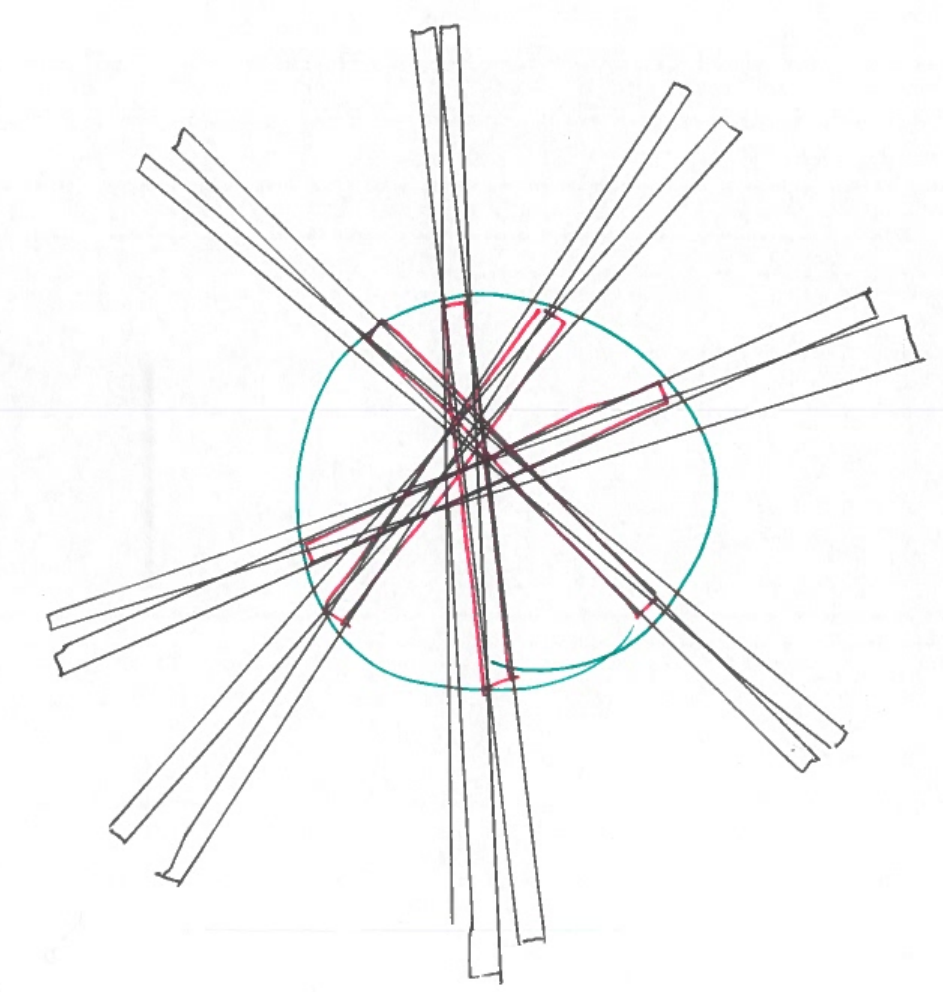}

In this picture,  the green circle is $B$,  the short red tubes belong to $\TT_B$,  and we see that each tube of $\TT_B$ comes from $\sim \mu(\TT_{T_\rho})$ longer tubes of $\TT$.   Therefore,  we can bound $\mu(\TT)$ by 

\begin{equation} \label{twoscales3}  \mu(\TT) \lessapprox \mu(\TT_{T_\rho}) \mu(\TT_B).  \end{equation}

We know that $\mu(\TT_{T_\rho}) \lessapprox | \TT_{T_\rho} |^\beta$,  and since we are in the not-sticky-at-all-scales case,  we also know that  $| \TT_{T_\rho}| \ll (\rho / \delta)^{2}$.  So we have

$$ \mu(\TT_{T_\rho}) \ll (\rho/\delta)^{2 \beta}. $$

Next, we have to bound $\mu(\TT_B)$.   Here it is much less clear what to do.   To get started,  let's consider the special case $\Delta_{max}(\TT_B) \lessapprox 1$.  

IF it happens that $\Delta_{max}(\TT_B) \lessapprox 1$, then we would get 

$$\mu(\TT_B) \lessapprox |\TT_B|^\beta \lessapprox (\rho^{-2})^\beta.$$

\noindent All together, we would have

$$ \mu(\TT)  \lessapprox \mu(\TT_{T_\rho}) \mu(\TT_B) \ll (\rho/\delta)^{2 \beta} \rho^{- 2 \beta} = \delta^{-2 \beta } \approx  |\TT|^\beta. $$

\noindent Since we assumed that $\mu(\TT) \approx | \TT |^\beta$,  this would give a contradiction.   

Now the hypothesis that $\TT_B$ has $\Delta_{max}(\TT_B) \lessapprox 1$ is a big IF (that's why I wrote IF in all caps).   The fact that $\Delta_{max}(\TT) \lessapprox 1$ does NOT imply that $\Delta_{max}(\TT_B) \lessapprox 1$.

\subsection{A surprising induction}

This is a key philosophical moment in the proof.   We are going to try to control $\TT_B$ using induction.   But the set of tubes $\TT_B$ need not obey $\Delta_{max}(\TT_B) \lessapprox 1$.   This makes it surprising to try to use $\TT_B$ in an inductive proof.   

To put this moment in context,  let us recall some background about the multilinear Kakeya problem.   
The multilinear Kakeya problem is a cousin of the Kakeya probelm that was solved by Bennett-Carbery-Tao in \cite{BCT}.   The proof was simplified in \cite{G},  and the proof there is only a few pages long.   Multilinear Kakeya involves $n$ sets of tubes $\TT_j$ in $\RR^n$,  where the tubes of $\TT_j$ are approximately parallel to the $x_j$ axis.   Multilinear Kakeya is important because it is much easier than Kakeya but still has many applications - for instance in work of Bourgain-Demeter on decoupling theory \cite{BD}. 

The key feature that makes multilinear Kakeya much easier than Kakeya is that if we intersect the tubes of each $\TT_j$ with a small ball $B$,  then the resulting sets of tubes $\TT_{j,B}$ obey the hypotheses of multilinear Kakeya.   Therefore,  we can easily apply induction to study the intersections of tubes in each small ball $B$.   In contrast,  the hypotheses of the Kakeya probem do not behave well when we restrict the tubes of $\TT$ to a small ball $B$.   

Nevertheless,  Wang and Zahl were able to use $\TT_B$ in an inductive argument.   While $\TT_B$ itself does not obey $\Delta_{max}(\TT_B) \lessapprox 1$,  Wang and Zahl manage to relate $\TT_B$ to some other sets of tubes $\tilde \TT$ with $\Delta_{max}(\tilde \TT) \lessapprox 1$,  and then we can use that $\mu(\tilde \TT) \lessapprox |\tilde \TT |^\beta$.   We will explore how they carry out this induction.

\vskip10pt

So far,  we have discussed the special case $\Delta_{max}(\TT_B) \lessapprox 1$.   Now we have to analyze the case $\Delta_{max}(\TT_B) \gg 1$.

If $\Delta_{max}(\TT_B) \gg 1$,  it means that for some convex set $K \subset B$,  $\Delta(\TT_B, K) \gg 1$.   Wang and Zahl consider the set $K$ that maximizes $\Delta(\TT_B, K)$.   The argument depends on the shape of $K$.   We will consider a few important shapes.

The first case we will consider is when $K = B$.   We will study this case in Section \ref{subsecdensity}.   Roughly we will show that in this case,  the Kakeya set $U(\TT)$ contains a large fraction of $B$,  and we will use this to close the induction.

The second case we will consider is the case when $K$ is a slab of dimensions $a \times r \times r$,  where $r$ is the radius of $B$ and $\delta \ll a \ll r$.   This case is the most important case in the argument.   We will discuss it in Section \ref{subsecslabcase}. 

We also briefly mention the case when $K$ is a very thin slab,  of dimensions $\approx \delta \times r \times r$.   Because of this case,  it helps to choose $\rho$ very close to $\delta$,  so that $r = \delta/\rho$ is very close to 1.   This case can now be handled by the $L^2$ method (without any induction).   
The intersection of tubes inside of each thin slab is basically the 2-dimensional Kakeya problem.  And the intersection of the thin slabs inside $B_r$ is well understood by the $L^2$ method.   Since $r$ is very close to 1,  this is sufficient to show that $|U(\TT)|$ is almost 1.

\subsection{Density of the Kakeya set in small balls} \label{subsecdensity}

Upper bounds for the multiplicity of the Kakeya set are equivalent to lower bounds for the volume of the Kakeya set.   In the course of the proof,  we will pay attention to the volume of the Kakeya set intersected with small balls. 

Let us write

$$U(\TT) = \cup_{T \in \TT} T.  $$

\noindent (Or if there is a shading $Y(T)$,  we would have $U(\TT) = \cup_{T \in \TT} Y(T)$.)

By double counting,  we have

$$ |U(\TT)| \approx \frac{ | \TT | |T| }{\mu(\TT)}.$$

Assuming $| \TT | \approx \delta^{-2}$,  

\begin{equation}
\mu(\TT) \lessapprox |\TT|^\beta \textrm{ is equivalent to } |U(\TT)| \gtrapprox \delta^{2 \beta}. 
\end{equation}

If $B_r$ is a ``typical' ball of radius $r$ intersecting the Kakeya set,  then we can write

$$ | U(\TT) | \approx | U(\TT_r) | \cdot \frac{ |U(\TT) \cap B_r| }{|B_r|}$$

We will refer to $ \frac{ |U(\TT) \cap B_r| }{|B_r|}$ as the density of the Kakeya set in $B_r$.  

Given that $\Delta_{max}(\TT) \lessapprox 1$ and $| \TT | \approx \delta^{-2}$,  it is not hard to show that for all $\delta \le r \le 1$ we have

\begin{equation}   | U(\TT_r) | \gtrapprox r^{2 \beta}. \end{equation}

Therefore,  in order to prove that $\mu(\TT) \ll | \TT |^\beta$,  it suffices to find a scale $r \gg \delta$ and prove a density lower bound

\begin{equation} \label{densitylb}
  \frac{ |U(\TT) \cap B_r| }{|B_r|} \gg  \left(\frac{\delta}{r} \right)^{2 \beta}.
  \end{equation}

Before we go on,  let us discuss the high density lemma in this context.   If $\TT$ is a Frostman set of $\delta$-tubes in $B_1$,  then it is not hard to check that $|U(\TT)| \gtrapprox \delta^{2 \beta}$.    The high density lemma,  Lemma \ref{lemmahighdensity},  gives an upper bound for $\mu(\TT)$ and hence a lower bound for $|U(\TT)|$.   If $\Delta_{max}(\TT) \gg 1$,  then it implies that $|U(\TT)| \gg \delta^{2 \beta}$.  

This argument allows us to handle the special case when $\Delta_{max}(\TT_B) \gg 1$ and $\Delta(\TT_B, K)$ is maximized when $K = B$.   In this case,  $\TT_B$ is Frostman in $B$.   Just using that $\TT_B$ is Frostman,  it follows that the density of $\TT_B$ in $B$ is $\gtrapprox (\delta/r)^{2 \beta}$.   But since $\Delta_{max}(\TT_B) \gg 1$,  the high density lemma implies that the density of $\TT_B$ in $B$ is $\gg  (\delta/r)^{2 \beta}$,  giving (\ref{densitylb}) and closing the induction in this case.

So far,  we have closed the induction in two special cases: the case when $\Delta_{max}(\TT_B) \lessapprox 1$ and the case when $\Delta_{max}(\TT_B) \approx \Delta(\TT_B, B) \gg 1$.

\subsection{When $\TT_B$ clusters into slabs} \label{subsecslabcase}

The most difficult and important case of the Wang-Zahl argument occurs when $\TT_B$ clusters into slabs.   This means that $\Delta_{max}(\TT_B) \gg 1$ and $\Delta(\TT_B, K)$ is maximized when $K$ is a slab of dimensions $a \times r \times r$,  where $r$ is the radius of $B$ and $a \ll r$.

We should mention that in the Heisenberg group example (over $\CC$),  the tubes of $\TT_B$ indeed cluster into slabs in this way.   Therefore,  this scenario was already regarded in the community as the most difficult and important scenario.   We will see that this scenario can also be handled by the ideas in the last couple subsections,  but with one new wrinkle.

\newcommand{\KK}{\mathbb{K}}

Let us first describe this scenario a little more carefully.   Suppose that $\KK$ is a set of $a \times r \times r$ slabs in $B$ so that $\Delta(\TT_B,  K) \approx \Delta_{max}(\TT_B)$ for each $K \in \KK$.   With some pigeonholing,  we can choose $\KK$ so that each $T_B \in \TT_B$ lies in $\sim 1$ slab $K \in \KK$.   Because the slabs of $\KK$ have maximal density,  it follows that $\Delta_{max}(\KK) \lessapprox 1$.   Then the $L^2$ method implies that $\mu(\KK) \lessapprox 1$.   Therefore,  we can study the slabs $K \in \KK$ one at a time.    For each $K \in \KK$,  define $\TT_K \subset \TT_B$ as $\TT_K = \{ T_B \in \TT_B: T_B \subset K \}$.   Now we have $\mu(\TT_B) \approx \mu(\TT_K)$.  

We know that $\Delta(\TT_K,  K) \gg 1$,  and we know that 

\begin{equation} \label{deffrostmanK}
\Delta(\TT_K,  K') \lessapprox \Delta(\TT_K,  K) \textrm{ for all convex } K' \subset K. 
\end{equation}

\noindent This setup is a natural generalization of the Frostman condition (\ref{Trhofrost}), so we say that $\TT_K$ is Frostman in $K$.

We might first try to study the density of small balls as in the last subsection.   The most direct version of this approach does not work,  but it can be made to work by adding an extra wrinkle.

The most natural radius to consider is $a$.   To close our induction,  it suffices to prove the density lower bound $\frac{ |U(\TT_K) \cap B_a|}{|B_a|} \gg \left(\frac{\delta}{a} \right)^{2 \beta}$.    If we average over all $B_a \subset K$,  it suffices to prove

\begin{equation} \label{goalhighdens}
\frac{ |U(\TT_K)|}{|K|} \gg \left(\frac{\delta}{a} \right)^{2 \beta}.
\end{equation}

At first,  we might hope to show that (\ref{goalhighdens}) follows from the fact that $\TT_K$ obeys the Frostman condition (\ref{deffrostmanK}) and that $\Delta(\TT_K, K) \gg 1$.   This does not work for the following reason.    Suppose that for every $\delta$,  there was a set $\TT_\delta$ of $\delta$-tubes in $B_1$ with $\Delta_{max}(\TT_\delta) \lessapprox 1$ and $|\TT_\delta| \approx \delta^{-2}$ and $\mu(\TT_\delta) \approx |\TT_\delta|^\beta$.    Then we could use these examples to build a bad set of tubes $\TT_K$ in $K$ so that

\begin{itemize}

\item $\Delta(\TT_K, K) \gg 1$.

\item $\TT_K$ obeys the Frostman condition (\ref{deffrostmanK})

\item $ \frac{ |U(\TT_K)|}{|K|} \approx \left(\frac{\delta}{a} \right)^{2 \beta}$

\end{itemize}

This bad example plays an important role in the story,  and we will take time to describe it.   Both the bad example and the way of dealing with it depend on understanding how tubes in $K$ relate to tubes in a unit ball.    

Recall that $K$ has dimensions $a \times r \times r$,  and each tube $T_K \in \TT_K$ has dimensions $\delta \times \delta \times r$.    Let $L$ be a linear map from the unit cube to $K$.   The map $L$ turns a tube of dimensions $s \times s \times 1$  into a plank of dimensions $as \times rs \times r$.    The most relevant scale is when the plank has thickness $\delta$,  so we set $a s = \delta$,  and so $s = \frac{\delta}{a}$.   Each $s$-tube in $B_1$ corresponds to a plank in $K$ of dimensions $\delta \times \frac{\delta}{a} r \times r$.    Not every plank in $K$ with these dimensions corresponds to an $s$-tube in $B_1$ -- the plank has to be oriented in an appropriate way.   Roughly speaking,  if a plank $P \subset K$ corresponds to a tube in $B_1$,  then the broadest face of $P$ should be approximately tangent to the broadest face of $K$.   We will call this class of planks horizontal planks.   To summarize we have a correspondence

$$ \textrm{ Horizontal planks in $K$ } \longleftrightarrow \textrm{ $s$-tubes in $B_1$}. $$

Here is a picture illustrating this correspondence:

\includegraphics[scale=.6]{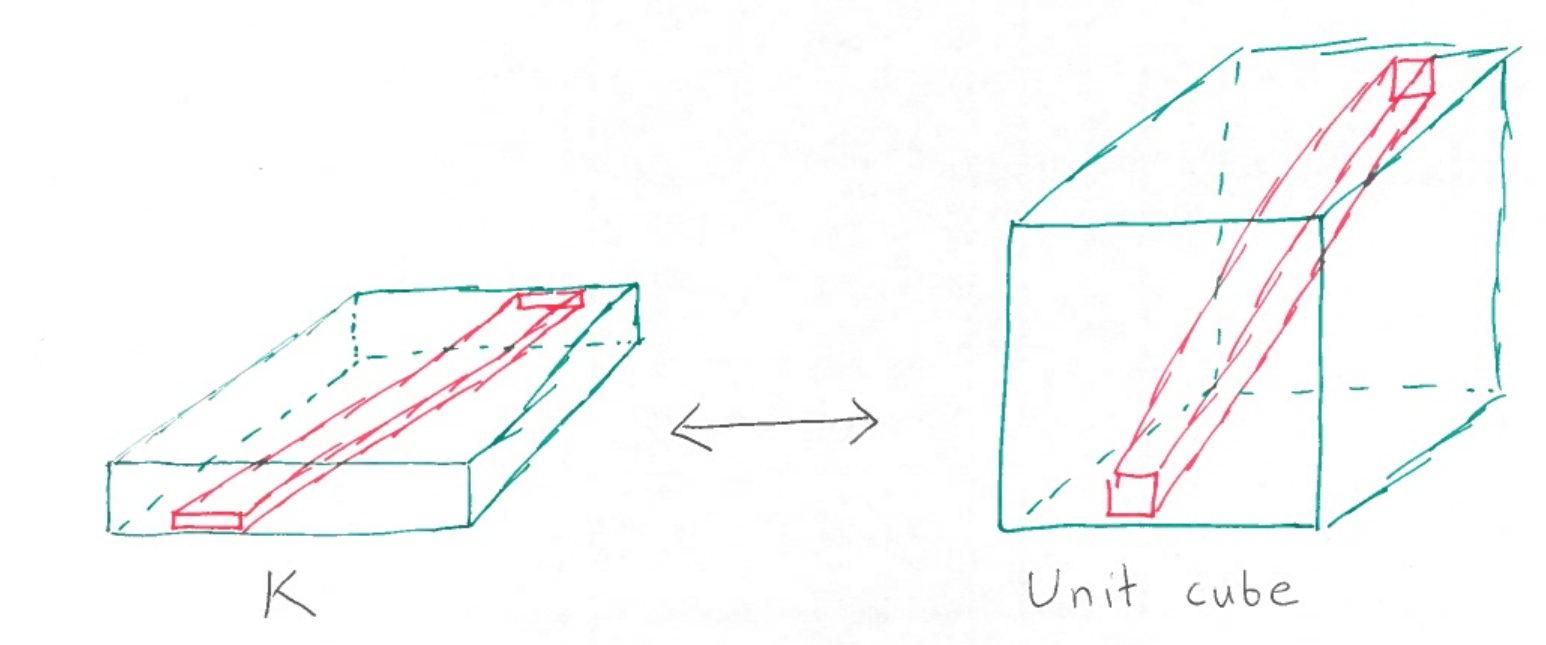}

The red plank on the left represents a horizontal plank $P$ in $K$.   Under the linear change of variables,  the plank $P$ corresponds to the red $s$-tube on the right.

Now we can describe our bad example.   Let $\TT_{s}$ with $\Delta_{max}(\TT_{s}) \lessapprox 1$ and $|\TT_{s}| \approx s^{-2}$ and $\mu(\TT_s) \approx |\TT_{s}|^\beta$.     Using our change of variables,  convert $\TT_{s}$ to a set of horizontal planks $\PP$ in $K$.   We have

$$ \frac{ |U(\PP)|}{|K|} = \frac{| U(\TT_{s}) |}{|B_1|} \approx s^{2 \beta} = \left( \frac{\delta}{a} \right)^{2 \beta}. $$

\noindent Since the dimensions of a plank $P$ are $\delta \times \frac{r}{a} \delta \times r$,  each plank $P \in \PP$ contains $\sim (r/a)^2$ distinct $\delta$-tubes.   Define $\TT_K$ to contain all the $\delta$-tubes in $P$ for each $P \in \PP$.    Now we have $\Delta(\TT_K,  K) \sim (r/a) \gg 1$.   It's not hard to check that $\TT_K$ is Frostman in $K$.   And $\frac{|U(\TT_K)|}{|K|} = \frac{ |U(\PP)|}{|K|} \approx (\delta/a)^{2 \beta}$.

This finishes the construction of our bad example.   The bad example shows that a direct generalization of the density argument from the last subsection does not work.

The bad feature of our bad example is that the tubes of $\TT_K$ cluster heavily in horizontal planks $P$.   The case when the tubes of $\TT_K$ cluster into horizontal planks is an important case in the Wang-Zahl argument.  In these notes, we will focus on this scenario.   
 We make this precise as follows.   We suppose $\PP$ is a set of horizontal planks.   For each $P \in \PP$,  let $\TT_{K,P} $ be a set of $\delta$-tubes in $P$.   The way that $\TT_{K,P}$ overlaps in $P$ is a 2-dimensional Kakeya type problem and the possible behaviors are well understood by the $L^2$ method.    We will consider the case that $U(\TT_{K,P}) = P$.  This case includes the bad example above,  in which $\TT_{K,P}$ consists of all the $\delta$-tubes in $P$.   

Wang and Zahl deal with this scenario by adding the following observation.  The angle between any two short tubes in the same horizontal plank $P$ is $\lesssim \delta/a$.   Therefore,  the way the tubes of $\TT_{K,P}$ intersect is related to the intersection of tubes in $\TT_{T_{\delta/a}}$,  and it can be controlled by studying $\TT_{T_{\delta/a}}$.   

For each plank $P$,  we let $\tilde P$ denote the corresponding $\delta/a$-tube in $B_1$.   And we let $\tilde \PP$ be the set of $\delta/a$-tubes corresponding to $P \in \PP$.   Since $\TT_K$ obeys the Frostman condition (\ref{deffrostmanK}),  it follows that $\PP$ is Frostman and so $\tilde \PP$ is also Frostman.    There are now two cases,  depending on $|\PP|$.   If $| \PP | \gg (\delta / a)^{-2}$,  then we are in the high density planks case.   If $| \PP | \lessapprox (\delta / a)^{-2}$,  then we are in the low density planks case.

\vskip10pt

{\bf High density planks case.}  In this case,  $\PP$ is Frostman in $K$ and $| \PP | \gg (\delta / a)^{-2}$.   Therefore,  $\tilde \PP$ is a Frostman set of $\delta/a$-tubes in $B_1$,  and $|\tilde \PP| \gg (\delta/a)^{-2}$.   In this case,  the high density lemma implies that $|U(\tilde \PP)| \gg (\delta/a)^{2 \beta}$,  and so

$$ \frac{ |U(\PP)|}{ |K|} = \frac{ |U(\tilde \PP)| }{|B_1|} \gg (\delta/a)^{2 \beta}. $$

\noindent We are focusing on the case that $U(\TT_{K,P}) = P$,  and so $U(\TT_K) = U(\PP)$.   Therefore,  we have $\frac{ |U(\TT_K)|}{ |K| } \gg (\delta/a)^{2 \beta}$,  which gives (\ref{goalhighdens}) and closes the induction.  

\vskip10pt

{\bf Low density planks case.}  In this case,  $\PP$ is Frostman in $K$ and $\Delta(\PP,  K) \lessapprox 1$,  and so $\Delta_{max}(\PP) \lessapprox 1$.  
In the low density planks case,  we will upper bound $\mu(\TT)$ using that

$$ \mu(\TT) \lessapprox \mu(\TT_{T_\rho}) \mu(\TT_B) \lessapprox \mu(\TT_{T_\rho}) \mu (\TT_{K,P}) \mu(\PP). $$

A typical point of our Kakeya set lies in $\mu(\TT_{K,P})$ short tubes from $\TT_{K,P}$.   The angle between two short tubes $T_K \in \TT_{K,P}$ is at most $\delta/a$.   Each of these short tubes belongs to $\mu(\TT_{T_\rho})$ long tubes.   Therefore,  $\mu(\TT_{T_\rho}) \mu(\TT_{K,P}) \lessapprox \mu(\TT_{T_{\delta/a}})$.  Plugging in,  we get

$$ \mu(\TT) \lessapprox \mu(\TT_{T_{\delta/a}}) \mu(\PP). $$

Now we bound $\mu(\PP) = \mu(\tilde \PP)$.  Since we are in the low density case,  $|\PP| = | \tilde \PP| \lessapprox (\delta/a)^{-2}$.   Since $\tilde \PP$ is Frostman,  we conclude that  $\Delta_{max}(\tilde \PP) \lessapprox 1$,  and so

$$ \mu(\PP) = \mu(\tilde \PP) \lessapprox | \tilde \PP|^\beta \lessapprox (\delta/a)^{- 2 \beta}. $$

So all together,  we have

$$ \mu(\TT) \lessapprox | \TT_{T_{\delta/a}}|^\beta (\delta/a)^{- 2 \beta}. $$

Since we are in the not-sticky-at-all-scales case,  $| \TT_{T_{\delta/a}} | \ll a^{-2}$,  and so, $\mu(\TT) \ll \delta^{-2 \beta} \approx | \TT|^\beta$. 
This closes the induction.

Notice that we really needed the assumptions of the not-sticky-at-all-scales case here.   We don't have any control of the scale $a$ except that $\delta \ll a \le r$.   The not-sticky-at-all-scales case tells us that $| \TT_{T_{\delta/a}}| \ll a^{-2}$ for all $a$ in this range.

\subsection{How do we reduce to the not-sticky-at-all-scales case?}

We saw in the argument above that it was important to reduce to the not-sticky-at-all-scales case.   In this subsection,  we explain how to do that.   We start by recalling the sticky case,  the not-sticky case,  and the not-sticky-at-all-scales case.

Recall that $\TT$ is a set of $\delta$-tubes in $B_1$ with $\Delta_{max}(\TT) \lessapprox 1$ and that $| \TT | \approx \delta^{-2}$. 

The sticky case means that for every $\rho \in [\delta,  1]$,  $\Delta_{max}(\TT_\rho) \lessapprox 1$.   Since $| \TT | \approx \delta^{-2}$,  this is equivalent to $| \TT_{T_\rho}| \approx (\delta/\rho)^{-2}$.

If $\TT$ is not sticky,  it means that there is some $\rho \in [\delta, 1]$ so that $| \TT_{T_\rho}| \ll (\delta/\rho)^{-2}$.   Such a $\rho$ must lie in the range $\delta \ll \rho \ll 1$.

We say that $\TT$ is not-sticky-at-all-scales if $| \TT_{T_\rho}| \ll (\delta/\rho)^{-2}$ for every $\rho$ in the range $\delta \ll \rho \ll 1$.

Here is the rough idea how to reduce the not-sticky case to the not-sticky-at-all-scales case.   Suppose that there is some scale $\rho$ so that $| \TT_{T_\rho}| \approx (\delta / \rho)^{-2}$.   It follows that $\Delta_{max}(\TT_\rho) \lessapprox 1$.   Now we can bound

$$ \mu(\TT) \lessapprox \mu(\TT_{T_\rho}) \mu(\TT_\rho), $$

\noindent which reduces our original problem to two similar problems at smaller scales.   We try to keep reducing in this way.   If one of the smaller problems is not-sticky-at-all-scales  then we are stuck and we cannot reduce further.   Otherwise we can reduce further.   If can keep reducing this way to very small problems,  it means that our original set of tubes $\TT$ was sticky,  and we can handle it using the sticky Kakeya theorem.   Otherwise,  we get stuck with a problem that is not-sticky-at-all-scales  and we can handle it using the argument we have sketched in this section.

(This reduction requires some care with small parameters to do in a rigorous way.)



\section{The high density lemma} \label{sechighdens}

The proof of the high density lemma, Lemma \ref{lemmahighdensity},  involves some important new ideas in its own right.   The high density lemma is about the sets of tubes that obey an important packing condition called the Frostman condition.   First let us recall this packing condition and put it in a more general context. 

Suppose that $\WW$ is a set of convex sets all contained in a given convex set $C$.   Recall that for a convex set $K$,  

$$\WW_K := \{ W \in \WW: W \subset K \}, $$

$$ \Delta(\WW, K) := \frac{ \sum_{W \in \WW_K} |W|}{|K|}. $$

We say that $\WW$ is Frostman in $C$ if $\WW$ is a set of convex sets in $C$ and 

\begin{equation} \label{deffrostman}
\Delta(\WW, K) \lessapprox \Delta(\WW,  C) \textrm{ for all } K \subset C
\end{equation}

Now we recall the statement of Lemma \ref{lemmahighdensity}:

\begin{lemma*}  (High density lemma) Suppose that $\beta$ is the critical exponent for sets of tubes with $\Delta_{max}(\TT)\lessapprox 1$,  defined in (\ref{critexp}).    Suppose that $\TT$ is a Frostman set of $\delta$-tubes in $B_1$.   Then 

$$ \mu(\TT) \lessapprox (\delta^{-2})^\beta (\delta^2 | \TT|)^{1 - \beta} . $$

\end{lemma*}

The exact bound on the right-hand side is not crucial.   The key point is to improve on the following trivial estimate.   Randomly decompose $\TT$ as a disjoint union $\TT = \cup_j \TT_j$,  where $| \TT_j | \sim \delta^{-2}$.   Since $\TT$ is Frostman,  it follows that for each $j$,  $\Delta_{max}(\TT_j) \lessapprox 1$,  and so $\mu(\TT) \lessapprox |\TT_j|^\beta = (\delta^{-2})^\beta$.   Now the trivial bound is 

\begin{equation} \label{trivbound}
\mu(\TT) \le \sum_j \mu(\TT_j) \lessapprox (\delta^{-2})^\beta (\delta^2 |\TT|)
\end{equation}

\noindent We see that if $\beta > 0$,  the bound in Lemma \ref{lemmahighdensity} improves on (\ref{trivbound}).  Even a small improvement is enough to power the inductive argument in the last section.   If (\ref{trivbound}) were sharp,  it would mean that for each $j$,  $|U(\TT)| \approx |U(\TT_j)|$.   This sounds intuitively unlikely: when $|\TT|$ is far bigger than $|\TT_j|$, we might expect $|U(\TT)|$ to be at least a little bigger than $|U(\TT_j)|$.   However,  it is not easy to prove this. 

Let $\gamma$ be the smallest exponent so that,  if $\TT$ is a Frostman set of tubes in $B_1$,  then

\begin{equation} \label{defgamma} \mu(\TT) \lessapprox (\delta^{-2})^\beta (\delta^2 | \TT|)^{\gamma}. \end{equation}

\noindent So Lemma \ref{lemmahighdensity} says that $\gamma = 1 - \beta$.   To power the inductive argument in the last section,  we just need to prove that if $\beta < 1$,  then $\gamma < 1 $.   We say that $\TT$ is a worse case example for Lemma \ref{lemmahighdensity} if

\begin{equation} \label{worsetcase} \mu(\TT) \approx (\delta^{-2})^\beta (\delta^2 | \TT|)^{\gamma}. \end{equation}

A crucial input to the proof is the sticky Kakeya theorem.   One scenario is that $\TT$ contains a sticky Kakeya set $\TT'$.   In this case,  $|U(\TT)| \ge |U(\TT')| \gtrapprox 1$,  and this leads to a very strong bound for $\mu(\TT)$.   When does $\TT$ contain a sticky Kakeya set?

Suppose that $1 = \rho_0 > \rho_1 > ... > \rho_N = \delta$ is a sequence of scales.   We can study $\TT$ at all these scales.   Recall that $\TT_\rho$ is the set of thickenings of $\TT$ at scale $\rho$.   And we let $ \TT_{\rho_j,  T_{\rho_{j-1}}}$ be the set of all $T_{\rho_j} \in \TT_{\rho_j}$ lying in the thicker tube $T_{\rho_{j-1}} \in \TT_{\rho_{j-1}}$.    If the sequence of scales $\rho_j$ is very dense,  and if each $\TT_{\rho_j,  T_{\rho_{j-1}}}$ is Frostman in $T_{\rho_{j-1}}$,  then $\TT$ contains a sticky Kakeya set.

This raises the question whether we can find such a sequence of scales.   We can try to build such a sequence of scales by adding one scale at a time,  and so it boils down to asking: is there a scale $\rho$ with $1 \gg \rho \gg \delta$ so that $\TT_\rho$ is Frostman and $\TT_{T_\rho}$ is Frostman?   Since $\TT$ is Frostman,  $\TT_\rho$ is automatically Frostman,  because whenever $K$ is a convex set that contains a $\rho$-tube,  $\Delta(\TT_\rho, K) \approx \Delta(\TT, K)$.   But $\TT_{T_\rho}$ are not necessarily Frostman.   We do know that for any convex set $K$,  $\Delta(\TT, K) \lessapprox \Delta(\TT,  B_1)$.   But if $\Delta(\TT,  T_\rho) \ll \Delta(\TT,  B_1)$,  then there could be some $K \subset T_\rho$ with $\Delta(\TT, K) \gg \Delta(\TT,  T_\rho)$.  

It is helpful to choose sets of maximal density.   For each $T_\rho$,  let $\WW(\TT_\rho)$ be a set of convex sets $W \subset \TT_\rho$ with maximal value of $\Delta(\TT,  W)$.    Because of the maximal density,  $\TT_W$ is Frostman in $W$ for each $W \in \WW(\TT_\rho)$.   After some pigeonholing and pruning the set $\WW(T_\rho)$, we can assume that each tube $T \in \TT_{T_\rho}$ lies in $\sim 1$ set $W \in \WW(\TT_\rho)$,  and each set $W \in \WW(T_\rho)$ contains about the same number of tubes $T \in \TT_{T_\rho}$.   

Now let $\WW = \bigcup_{T_\rho \in \TT_\rho} \WW(T_\rho)$.   Each $W \in \WW$ has dimensions $a \times b \times 1$,  where $\delta \le a \le b \le \rho$.   Note that if $W \in \WW(T_\rho)$,  then $T_\rho$ is essentially the $\rho$-neighborhood of $W$,  and so each $W$ belongs to $\WW(T_\rho)$ for $\sim 1$ $T_\rho \in \TT_\rho$.   So each tube $T \in \TT$ lies in $\sim 1$ set $W \in \WW$.   And by pigeonholing,  we can assume that each $W \in \WW$ contains around the same number of tubes $T \in \TT$.

Since $\TT$ is Frostman in $B_1$,  it follows that $\WW$ is also Frostman in $B_1$.   

Now we can bound $\mu(\TT)$ by

\begin{equation} \label{twoscales4} \mu(\TT) \lessapprox \mu(\TT_W) \mu(\WW).
\end{equation}

\noindent Since $\TT_W$ and $\WW$ are both Frostman,  we can try to bound each factor by induction.

We note that the sets $W$ are convex sets but not necessarily tubes.    After some pigeonholing,  we can assume that each $W \in \WW$ has dimensions $\sim a \times b \times 1$ with $\delta \le a \le b \le 1$. 
Therefore,  we have to generalize the setup of our problem to include not just intersecting tubes but intersecting convex sets.   Wang and Zahl systematically develop a theory of intersecting convex sets instead of just intersecting tubes.  One extreme example is when each $W$ is a slab of dimensions $a \times 1 \times 1$.   Sharp bounds for intersecting slabs have been known for a long time by the $L^2$ method.   So the case of slabs is fully understood.   If $W$ has dimensions $a \times b \times 1$ with $a \ll b \ll 1$,  then the shape of $W$ is intermediate between a tube and a slab.   Leveraging the sharp bounds for slabs,  the bounds for this intermediate case are stronger than the bounds for actual tubes.   Therefore,  if $a \ll b$,  then when we apply induction to (\ref{twoscales4}),  we get a gain.   If we assume that $\TT$ was a worst case example for Lemma \ref{lemmahighdensity},  then this scenario cannot occur and we must essentially have $a = b$.

Another scenario is that $W$ could be an original tube $\delta \times \delta \times 1$.   In this case,  we have not decomposed our original problem at all.   If we choose $\rho$ so that $\Delta(\TT,  T_\rho) \gg 1$,  then we get $\Delta(\TT, W) \gtrsim \Delta(\TT,  T_\rho) \gg 1$,  and so $W$ must be significantly larger than a single tube.

To summarize,  if $\TT$ is a worse case example for the high density lemma,  and if $|\TT| \gg \delta^{-2}$,  then  there is a scale $a$ with $\delta \ll a \ll 1$ so that $\TT_{T_a}$ is Frostman and $\TT_a$ is Frostman.   Moreover,  $\TT_a$ and $\TT_{T_a}$ will also be worse case examples for the high density lemma,  which helps to iterate.

On the other hand,  if $\TT$ is a Frostman set of tubes in $B_1$ with $|\TT| \approx \delta^{-2}$,  then we have $\Delta_{max}(\TT) \approx 1$ and so we can  bound $\mu(\TT) \lessapprox (\delta^{-2})^\beta$.

Iterating this argument,  we can choose a sequence of scales $1 = \rho_0 \gg \rho_1 \gg \rho_2 \gg ... \gg \rho_N = \delta$ so that  $ \TT_{\rho_j,  T_{\rho_{j-1}}}$ is always Frostman, and for each $j$ one of the following holds:

\begin{enumerate}

\item $\rho_j / \rho_{j+1}$ is very small, or 
\item  $ \Delta_{max}( \TT_{\rho_j,  T_{\rho_{j-1}}}) \lessapprox 1$, and so $\mu(  \TT_{\rho_j,  T_{\rho_{j-1}}}) \lessapprox (\rho_{j-1}/ \rho_j)^{2 \beta}. $

\end{enumerate}

Here is a picture showing how this sequence of scales may look:

\includegraphics[scale=.7]{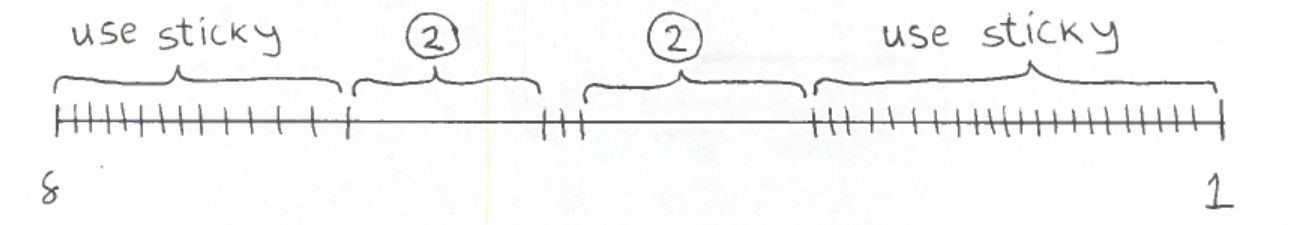}

Here each short vertical line represents a scale $\rho_j$.   These scales are generally quite close together,  except for two significant gaps.  Each significant gap must be in case (2),  and so we labelled them (2).  

If $|\TT| \gg \delta^{-2}$,  not every interval can be in case 2 --  a definite fraction of intervals must be in case 1.   Since the intervals in case 1 are very small,  a definite fraction of scales must consist of many small intervals.   On any such block of small intervals,  we can apply sticky Kakeya,  giving a very strong bound.   In our picture,  we have drawn two such blocks of small intervals,  and they are labelled ``use sticky Kakeya''.

To bound $\mu(\TT)$,  we begin by factoring $\mu(\TT)$ into contributions coming from different scale ranges.   For instance,  in the scenario illustrated above,  we would first name key scales as in the picture below:

\includegraphics[scale=.7]{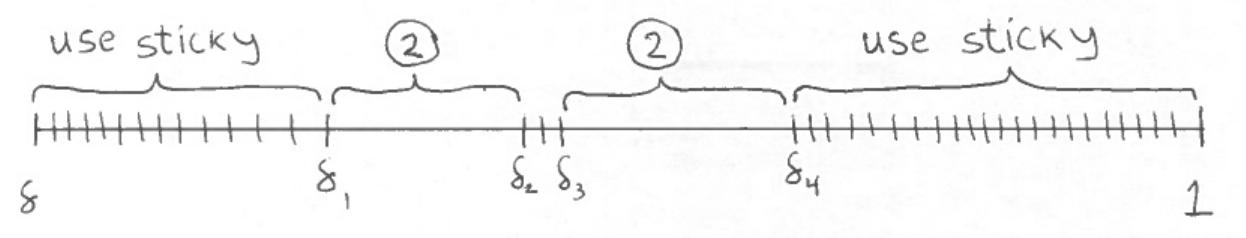}

\noindent Now we can bound $\mu(\TT)$ by

\begin{equation} \label{multscalefactor} \mu(\TT) \lessapprox \mu( \TT_{T_{\delta_1}}) \mu(\TT_{\delta_1,  T_{\delta_2}})  \mu(\TT_{\delta_2,  T_{\delta_3}})  \mu(\TT_{\delta_3,  T_{\delta_4}})  \mu(\TT_{\delta_4}) . \end{equation}

\noindent The five factors on the right-hand side correspond to five scale ranges in the picture.    Note that if we used the trivial bound (\ref{trivbound}) to bound each factor on the right-hand side,  then we would get back the trivial bound.   Instead,  we use sticky Kakeya on the scale ranges labelled sticky,  and we use item (2) above on the scale ranges labelled (2).   This leaves some small scale ranges  where we have to use the trivial bound (\ref{trivbound}),  corresponding to $\gamma = 1$.   Since $\beta < 1$,  the bound in the sticky case is better than the trivial bound,  and so our overall bound for $\mu(\TT)$ is better than the trivial bound (\ref{trivbound}).   A careful calculation gives the bound $\gamma = 1 - \beta$.  

This finishes our outline of the proof that sticky Kakeya implies Kakea.  

\vskip10pt

Notice that in this argument,  we broke the scales from $\delta$ to 1 into several ranges in a strategic way.   The idea of choosing these ranges strategically was introduced by Keleti and Shmerkin in \cite{KS},  and it has become a major tool in this circle of problem.   For instance,  it plays a key role in the solution of the Furstenberg set problem in \cite{KS} and \cite{RW}.

\vskip10pt

Over the last two sections, we have outlined the argument from \cite{WZ},  reducing the general case of Kakeya to the sticky case.   It can be considered as an argument by induction on scales.   We define $\beta$ as the best exponent so that if $\Delta_{max}(\TT) \lessapprox 1$,  then $\mu(\TT) \lessapprox |\TT|^\beta$.    If $\TT$ is not sticky,  we show (using induction),  that $\mu(\TT) \ll | \TT |^\beta$.    This shows that the worst case is sticky,  and so reduces the general Kakeya problem to sticky Kakeya.  

The set of tubes $\TT_B$ played a key role in this inductive argument.   Recall that we defined $\TT_B$ by intersecting the tubes of $\TT$ with a smaller ball $B$.  Since we don't know whether $\Delta_{max}(\TT_B) \lessapprox 1$,  it sounded difficult to involve $\TT_B$ in an inductive argument.  

How did we control $\TT_B$ using induction?   In the argument,  $\TT_B$ is treated as an arbitrary set of tubes.  There are several cases.   In the most interesting case,  we related $\TT_B$ to a set of planks $\PP$.   After a linear change of variables, $\PP$ becomes a set of tubes $\TT'$ obeying a Frostman condition,  and we estimate $\mu(\TT')$ using the high density lemma.   When we look inside the proof of the high density lemma,  we examine $\tilde \TT'$ at various scales.  We eventually find some sets of tubes $\tilde \TT$ with $\Delta_{max}(\tilde \TT) \lessapprox 1$.   We can apply our inductive assumption to each such $\tilde \TT$.   In this subtle way,  we were able to study $\TT_B$ using induction.

\section{What about the Katz-Zahl example?}

As we mentioned in Section \ref{secstickynot},  Katz and Zahl found a cousin problem to the Kakeya problem where the analogue of Kakeya does not hold but the analogue of the sticky case appears likely to hold.   That example made me think it was unlikely that the Kakeya problem could be reduced to the sticky case.   Wang and Zahl did reduce the general Kakeya problem to the sticky case,  and so it is natural to ask why their method does not apply to the Katz-Zahl example.

The Katz-Zahl example concerns a cousin of the Kakeya problem where $\RR$ is replaced by the ring $ A= \FF_p[x] / (x^2)$.   The ring $A$ has a natural notion of distance with two distinct length scales.   If $a + bx \in A$ with $a,b \in \FF_p$,  we define

$$ \| a + b x \|_A := \begin{cases} 1 & \textrm{ if } a \not= 0 \\ p^{-1} &  \textrm{ if } a = 0,  b \not= 0 \\ 0  & \textrm{ if } a = b = 0 \\ \end{cases} $$

\noindent There is a cousin of the Heisenberg group in $A^3$ and it leads to a counterexample to the analogue of Theorem \ref{main}.   But unlike in $\CC^3$,  the Heisenberg group cousin in $A^3$ is {\it not} sticky.   It appears likely that in $A^3$,  the sticky case of Wolff axiom Kakeya conjecture is true,  but the general conjecture is false.

Looking back at the proof of Wang-Zahl,  the key distinction between $\RR$ and the ring $A$ is that the ring $A$ has only two distinct non-zero scales.   The proof of the Kakeya conjecture requires discussing many scales in order to run the multiscale analysis.   

(When we run the Wang-Zahl argument on the Katz-Zahl example,  I think that it shows that the intersection of $U(\TT)$ with a ball of radius $p^{-1}$ consists of a union of slabs of dimensions $\delta \times p^{-1} \times p^{-1}$.    The Katz-Zahl example has this property.   When we study Kakeya over $\RR$,  we would also consider balls of other radii, like maybe $p^{-1/100}$.  And in some cases of the argument we would show that the intersection of $U(\TT)$ with a ball of radius $p^{-1/100}$ consists of a union of slabs of dimensions $\delta \times p^{-1/100} \times p^{-1/100}$.   However,  working over the ring $A$,  we cannot carry this out because there is no such scale $p^{-1/100}$.)

There is an analogous ring $A_N$ which has many scales.   Define $ A_N = \FF_p[x] / (x^N)$.   An element $r$ in the ring $A_N$ can be written as $r = \sum_{j=0}^{N-1} a_j x^j$ with $a_j \in \FF_p$.   If $j(r)$ is the least $j$ so that $a_j \not= 0$,  then $\| r \|_{A_N} = p^{-j}$.   I think that if $N$ is sufficiently large,  then the argument reducing Kakeya to sticky Kakeya has a good chance to work in the ring $A_N$.

More broadly,  we can ask: ``For which problems of Kakeya type is the worst case sticky?''  From the Katz-Zahl example,  we learn that we should consider problems with many different scales.   As far as I know,  it looks plausible that for a broad range of Kakeya type problems involving many scales,  the worst case is sticky.    On the other hand,  in the Furstenberg set problem,  the known worst case example comes from a grid,  and it is not sticky.   It looks plausible that examples which are sticky in a strong sense cannot match the worse case in the Furstenberg set problem.   
It would be interesting to explore how much the proof from \cite{WZ} can be generalized.

\section{What about higher dimensions?}

It is natural to ask about generalizing this work to higher dimensions.   No one knows yet how to do that.   There is one key difference in dimension $n \ge 4$.   Our whole discussion has been based around the convex Wolff axioms.   But it is well known that in dimension at least 4,  the analogue of the Kakeya conjecture with convex Wolff axioms is not true.   Let us write down a precise question and then recall the example.

Our definitions of density and $\Delta_{max}$ make sense in any dimension.

Suppose that $\TT$ is a set of $\delta$ tubes in $\RR^n$.   For any convex set $K \subset \RR^n$,  we define the density the density of $\TT$ in $K$ as

$$ \Delta(\TT, K) = \frac{ \sum_{T \in \TT_K} |T|}{ |K| }. $$

\noindent Then we define

$$ \Delta_{max}(\TT) := \max_K \Delta(\TT, K). $$

\noindent And we define the typical multiplicity of $\TT$ as

$$ \mu(\TT) = \frac{ \sum_{T \in \TT} |T|}{ | \cup_{T \in \TT} |}. $$

In dimension 4,  there is an example of a set of $\delta$-tubes $\TT$ with $\Delta_{max}(\TT) \lesssim 1$ and $\mu(\TT) = \delta^{-1}$.   In this example,  the tubes are all contained in the $\delta$-neighborhood of a degree 2 algebraic hypersurface,  such as

$$Z = \{ x \in \RR^4:  x_1^2 + x_2^2 - x_3^2 - x_4^2 = 1 \} $$

\noindent A typical point of $Z$ lies in a 1-parameter family of lines in $Z$.   For instance,  if we intersect $Z$ with the hyperplane $x_1 = 1$,  then we get a 2-dimensional cone defined by $x_2^2 - x_3^2 - x_4^2 = 0$,  and this cone contains a 1-parameter family of lines.   So the point $(1,0,0,0)$ lies in a 1-parameter family of lines in $Z$.   But the variety $Z$ is very symmetric.   The Lie group $O(2,2)$ acts transitively on $Z$.   Since this is a linear action,  the action maps lines to lines.   And so every point of $Z$ lies in a 1-parameter family of lines in $Z$.   

The tubes of $\TT$ are $\delta$-neighborhoods of these lines.   There are $\delta^{-3}$ tubes of $\TT$.   The set $U(\TT)$ is the $\delta$-neighborhood of $Z$,  and so $|U(\TT)| \sim \delta$ and so $\mu(\TT) \sim \delta^{-1}$.  

The tubes of $\TT$ cluster near the algebraic surface $Z$,  but they don't cluster in any convex set.   It is not difficult to check that for any convex set $K$,  $\Delta(\TT, K) \lesssim 1$,  and so $\Delta_{max}(\TT) \lesssim 1$.

There are two possible ways around this example.   One way is to focus on the original version of the Kakeya conjecture with direction-separated tubes.   In this setup,  we suppose that $\TT$ is a set of $\delta^{-(n-1)}$ $\delta$-tubes in $\RR^n$ in $\delta$-separated directions.   (The example above is not a counterexample to this direction separated Kakeya conjecture,  because the directions of the tubes in the example all lie near a degree 2 algebraic surface in $S^3$. )   Working with convex Wolff axioms played an important role in the reduction to the sticky case.   Can it be made to work in the directon-separated framework?

The second way around this example is to replace the convex Wolff axioms with more general polynomial Wolff axioms that take account of how tubes of $\TT$ pack into neighborhoods of low degree algebraic varieties.  I think this is a natural direction to pursue.   But replacing convex sets by semi-algebraic sets is a significant issue and is not just a technical modification of the 3-dimensional proof.

\section{Next steps}

As I mentioned at the start of this survey,  I am currently working on digesting and checking the proof.    My next step is to write a detailed outline of the proof in \cite{WZ},  which includes all the cases and all the steps.   I hope to have something that I can share by end of June.

\end{document}